\documentclass[12pt,twoside]{article}
\usepackage[english]{babel}
\usepackage[latin1]{inputenc}
\usepackage[T1]{fontenc}

\usepackage[pdftex]{graphicx}
\usepackage{graphics}

\usepackage{amsmath}
\usepackage{amssymb}
\usepackage{amsthm}
\usepackage{amsfonts}
\usepackage{amssymb}
\usepackage{amsmath}
\usepackage[mathscr]{eucal}
\usepackage{algorithm}
\usepackage{algorithmic}
\usepackage{amscd}
\usepackage{amsthm}

\usepackage{comment}

\usepackage{amsmath,amsthm}

\usepackage{array}
\usepackage{color}

\DeclareMathOperator*{\argmin}{arg\,min}

\theoremstyle{plain}
\newtheorem{theorem}{Theorem}[section]
\newtheorem{Lemma}[theorem]{Lemma}
\newtheorem{proposition}[theorem]{Proposition}
\newtheorem{Corollary}[theorem]{Corollary}

\theoremstyle{definition}
\newtheorem{definition}{Definition}[section]

\theoremstyle{remark}

\usepackage{url}
 \pdfoutput=1\pdfcompresslevel=9\pdfadjustspacing=1
 \usepackage[pdftex,bookmarks,a4paper,colorlinks]{hyperref}
 \usepackage[pdftex]{thumbpdf}

 \hypersetup{bookmarksnumbered,plainpages=false,hypertexnames=false}
 \hypersetup{pagecolor=blue,linkcolor=blue,citecolor=blue,urlcolor=blue}
  \hypersetup{pdfcreator=PDFLaTeX + Hyperref}

\title{A data driven trimming procedure for  robust classification}

\begin{document}

\author{Marina Agull\'o Antol\'in$^*$, Eustasio del Barrio\footnote{{IMUVa. Universidad de Valladolid. 7, paseo de Bel\'en, 47011 Valladolid. Spain. These authors have been partially supported by the Spanish Ministerio de Econom\'ia y
Competitividad, grants MTM2014-56235-C2-1-P, and MTM2014-56235-C2-2,
and by Consejer\'ia de Educaci\'on de la Junta de Castilla y Le\'on,
grant VA212U13.}} and Jean-Michel Loubes\footnote{Institut de Math\'ematiques de Toulouse, Universite Paul Sabatier. 118, route de Narbonne F-31062 Toulouse Cedex 9, loubes@math.univ-toulouse.fr}}
%
%

\maketitle


\begin{abstract}
Classification rules can be severely affected by the presence of
disturbing observations in the training sample. Looking for an
optimal classifier with such data may  lead to unnecessarily complex
rules. So, simpler effective classification rules could be achieved
if we relax the goal of fitting a good rule for the whole training
sample but only consider  a fraction of the data. In this paper we
introduce a new method based on trimming to produce classification
rules with guaranteed performance on a significant fraction of the
data. In particular, we provide an automatic way of determining the
right trimming proportion and obtain in this setting oracle bounds
for the classification error on the new data set.
\end{abstract}


\section{Introduction}
In the usual classification setting we observe a collection of pairs
of i.i.d copies $(Y_i,X_i) \in \{0,1\} \times \mathbb{R}^p$ with
$i=1,\dots,n$ of a random variable $(Y,X)$ with distribution $P$.
$Y$ is the label to be forecast according to the value of the
variables $X$. A classifier is a function $g: \mathbb{R}^p \mapsto
\{0,1\}$ that predicts the label of an observation. An observation
is misclassified if $Y \neq g(X)$. Hence, the performance of a
classifier can be measured by its classification error defined as
$R(g)=P((y,x)\in\{0,1\} \times \mathbb{R}^p:y\neq g(x))$. During the
last decades, the classification problem a.k.a pattern recognition
has been extensively studied and there exists a large variety of
methods to find optimal classifiers in different settings. We refer
for instance to \cite{ElementsStatisticalLearn},
\cite{SVMCristianini} or \cite{LugosiPattern} and references therein
for a survey.

When the number of observations grows large or in a high dimensional
case, some of the data may contain observation errors and may be
considered as contaminating data. The presence of such observations,
if not removed, hampers the efficiency of classifiers since many
classification methods are very sensitive to outliers. Actually, if
the learning set is too corrupted,  training a classifier over this
set leads to bad classification rates. Hence there is a growing need
for robust methods to tackle such issue.  Pioneered in
\cite{MR0161415}, we refer to  \cite{MR829458} for a review of
robust methods.

A solution to cope with this issue is to allow the classifier not to
label all points but to reject some observations that may seem too
difficult to be classified. This point of view is studied in
\cite{herbei2006classification} and
\cite{bartlett2008classification}. Another  general idea  is to
remove a proportion of contaminating data to guarantee the
robustness of the method. Such data are defined as outliers in the
sense that they are far from the model used to generate the data.
Yet detecting automatically outliers is a difficult task since their
mere definition is unclear and highly depends on each particular
case. Much research has been done in this direction and many
analysis provide several ways of determining whether an observation
is an outlier. For instance in \cite{MR2752147}, in the case of SVM
classification, the author proposes to remove observations using an
outlier map. In \cite {MR2248007}, the authors rely on a function
that measures the impact of contamination of the distribution on the
classification error obtained by minimizing the empirical risk
criterion. In a regression framework, Lasso estimators suffer from
outliers. They can also be modified in order to enhance robustness
as in \cite{MR2808669}, \cite{maronna2011robust} or in \cite{Alfons}
where the authors discard the points for which the residuals are the
largest.

In a probabilistic framework, removing observations that achieve bad
classification error, corresponds to trimming the initial
distribution of the observations and replacing it by a similar
distribution $Q$ up to some data that will be considered as outliers
for the classification rule. Trimming methods for data analysis have
been described firstly in \cite{rousseeuw1984least} and later some
statistical properties are studied in \cite{MR1983163} or in
\cite{MR1439314}. Yet very few theoretical results exist to study
how to choose the actual boundary between an acceptable observation
and an outlier.  Moreover, in that case little is known about
whether this choice modifies the classification error or how. For
both theoretical and practical purposes, such a choice must be
guided in order to take into account the amount of variability
generated by the corrupted data.

In this paper we provide some theoretical guarantee to choose the
level of data to be removed. For this we consider the set of trimmed
distributions obtained from the initial distribution of the data and
look for an automatic rule that reduces the classification error of
a collection of classifiers by removing some properly selected
observations. The more data is removed, the easier it becomes to
classify the data, leading to a perfect classification if the
classification rate is small enough. Yet, removing too many data
reduces the interest of the classification procedure. If too many
observations are left aside then the chosen classification rule may
be good for distribution which is possibly very far from the true
distribution of the data. We provide in this work an empirical rule
that automatically selects the minimum level of trimming to reduce
the classification error for a class of classifiers. Simultaneously,
the best classifier for the trimmed set of observations is chosen
among a collection of classification rules and for this, we prove an
oracle inequality that governs the statistical properties of this
methodology. \vskip .1in

The paper falls into the following parts. Section~\ref{s:pct} is
devoted to the description of the probabilistic framework of outlier
selection using trimming distributions. We precisely define the
trimmed classification errors and their relationships with the usual
classification errors for both the empirical and theoretical error.
In section~\ref{s:choice} we provide the automatic selection rule
for the trimming level and the best trimmed classifier for which we
provide and  oracle inequality. This model selection result is
illustrated  with the case of linear classifiers.
Section~\ref{s:conclu} provides some conclusions and perspectives
for these results. The proofs and some technical results are
gathered in the Appendix.

\section{Partial Classification with trimming} \label{s:pct}

As in the introduction, we assume that we observe an $n$ i.i.d
sample $(Y_i,X_i)_{i=1,\dots,n} \in \{0,1\} \times \mathbb{R}^p$
with distribution $P$. Set $g: \mathbb{R}^p \mapsto \{0,1\}$, a
classification rule, we denote the classification error as
$$R(g)=P((y,x)\in\{0,1\} \times \mathbb{R}^p:y\neq g(x)).$$
Since the underlying distribution of the observations is unknown,
the classification error  $R$ is estimated by its empirical
counterpart, the empirical error defined as
\begin{equation*}\label{DefErrEmp}
R_n(g):=\frac{1}{n}\displaystyle\sum_{i=1}^nI_{(g(X_i)\neq Y_i)},
\end{equation*}
where $I_{(g(X)\neq Y)} = 1$ if $g(X)\neq Y$ and $0$ otherwise. \\
\indent Trimming  a data sample of size $n$ is usually defined as
discarding a given fraction of the data while reweigthing the other
part. Let $\alpha$ be the proportion of observations we can trim,
and consider that $n\alpha=k\in\mathbb{N}$. Then, trimming consists
of removing $k$ observations and giving weight $1/(n-k)$ to the
rest. Among all the possible trimmings, we will call \emph{empirical
trimmed classification error} the one that minimizes the sum

\begin{equation}\label{DefErrorEmpiricoRecortado}
R_{n,\alpha}(g):=\min_{w \in W} \displaystyle\sum_{j=1}^n w_i I_{(g(x_j)\neq
y_j)}
\end{equation}
with
\begin{equation}\label{WeightsConditions}
W=\{w=(w_1,\ldots,w_n)/ \ 0\leq w_i \leq \frac{1}{n(1-\alpha)}; \
i=1,\ldots,n \ \wedge \ \displaystyle\sum_{i=1}^n w_i=1\}.
\end{equation}
To study the theoretical counterpart of this quantity, which we will
call trimmed classification error, we will consider the set of
trimmed distributions as follows. From a probabilistic point of
view, trimming a distribution consists of replacing the initial
distribution of the observations by a new measure built by a partial
removal of points in the support of the initial distribution. We
thus can provide the following definition for the the trimming of a
distribution. Here, $\mathcal{P}$ denotes the set of probabilities
on $\{0,1\} \times \mathbb{R}^p$.

\begin{definition}
Given $\alpha\in(0,1)$, we define the set of $\alpha$-trimmed
versions of P by
\begin{equation*}\label{TrimmDef}
    \mathcal{R}_\alpha(P):=\left\{Q\in \mathcal{P} :\, Q\ll P,\, \frac{dQ}{dP}\leq\frac{1}{1-\alpha}
    \ P-a.s.\right\}.
\end{equation*}
\end{definition}
This entails that a trimmed distribution $Q \in
\mathcal{R}_\alpha(P)$ can be seen as a close modification of a
distribution $P$ obtained by removing a certain quantity of data
(see \cite{alvarez2012similarity} and \cite{alvarez2012trimmed}).
When dealing with a classification rule, one is interested in
looking for the data for which the classification rule performs
well. Hence, we aim at improving the classification error by
changing the underlying distribution of the observations using a
trimming scheme in order to modify the distribution but yet in a
controlled, limited way. With this goal we introduce the trimmed
classification error for a rule $g: \mathbb{R}^p \mapsto \{0,1\}$.

\begin{definition}
Given $\alpha\in(0,1)$, we define the \emph{trimmed classification
error} of a rule as the infimum of the $\alpha$-trimmed
probabilities of misclassifying future observations
\begin{equation*}\label{ClassTrimError}
    R_\alpha(g):=\displaystyle\inf_{Q\in\mathcal{R}_\alpha(P)}Q(g(x)\neq
    y).
\end{equation*}
\end{definition}

There is a simple relation between the trimmed classification error
and the general classification error as the next result shows.

\begin{proposition}\label{R_alpha(g)}
Given a trimming level $\alpha\in(0,1)$ and a classification rule
$g$,
\begin{equation}\label{RalphaFunR}
    R_\alpha(g)=\frac{1}{1-\alpha}\left(R(g)-\alpha\right)_+.
\end{equation}

\end{proposition}
This proposition shows the effect of trimming on classification. We
write $g_B$ for the Bayes classifier, namely, the classification
rule that yields the minimal classification error. We also write
\begin{equation*}
    Err(P):=\displaystyle\min_g R(g)=R(g_B),
\end{equation*}
for the \emph{Bayes classification error}. The optimal trimming for
classification removes the misclassified points in such a way that
if the classification error is less than the percentage of points
that can be removed, then all the points are classified without
error.

Similar to the Bayes rule, we can define a trimmed Bayes
classification rule and the trimmed Bayes error as follows.

\begin{definition}
An \emph{$\alpha$-trimmed Bayes classifier} or
\emph{$\alpha$-trimmed Bayes classification rule} is a classifier
that achieves the minimum $\alpha$-Trimmed classification error
\begin{equation*}
    g_B^\alpha:=\displaystyle\argmin_g R_\alpha(g).
\end{equation*}
The corresponding classification error is thus the
\emph{$\alpha$-Trimmed Bayes error} defined as
\begin{equation*}\label{BayesRecortado}
    Err_\alpha(P):=\displaystyle\inf_{Q\in\mathcal{R}_\alpha(P)}Err(Q)=\displaystyle\min_g R_\alpha(g)=R_\alpha(g^{\alpha}_B).
\end{equation*}
\end{definition}
The following proposition compares these two errors.
\begin{proposition} \label{RelBayesrec}
\begin{equation*}\label{th:general}
\mbox{Err}_\alpha(P)=\frac{(R(g_B)-\alpha)_+}{1-\alpha}=\frac{(\mbox{Err}(P)-\alpha)_+}{1-\alpha}.
\end{equation*}
\end{proposition}

If $\mbox{Err}(P)\leq \alpha$ then $\mbox{Err}_\alpha(P)=0$, but if
$\mbox{Err}(P)> \alpha$ then
$\mbox{Err}_\alpha(P)=\frac{\left(\mbox{Err}(P)-\alpha\right)_+}{1-\alpha}>0$,
which indicates that $\mbox{Err}_\alpha(P)=0$ is equivalent to $
\mbox{Err}(P)\leq \alpha.$ This means, the minimum $\alpha$ which
gives us the perfect separation is the value that corresponds to the
Bayes error.

Usually we do not look for the optimum classifier among all possible
classifiers, but we restrict ourselves to a smaller class of
classifiers. Let $\mathcal{F}$ be a class of classifiers and let
${f^\star}\in\mathcal{F}$ be the classifier which gives us the
minimum classification error within the class. We denote as
$R(\mathcal{F})$ the minimum classification error in $\mathcal{F}$,
that is
$$R(\mathcal{F}):=\displaystyle\min_{f\in\mathcal{F}}R(f)=R({f^\star}).$$
In the same way we denote the trimmed error of the class
$\mathcal{F}$ as $R_\alpha(\mathcal{F})$. Hence, given Proposition
\ref{R_alpha(g)},
$$R_\alpha(\mathcal{F}):=\displaystyle\min_{f\in\mathcal{F}}R_\alpha(f)=\displaystyle\min_{f\in\mathcal{F}}\frac{\left(R(f)-\alpha\right)_+}{1-\alpha}.$$
The classifier that minimizes $R(f)$ also minimizes this quantity
and so the classifier that minimizes the error in the class
$\mathcal{F}$ is a minimizer of the trimmed error in the class.

Proposition \ref{R_alpha(g)} can be trivially applied to the
empirical trimmed classification error introduced in
\eqref{DefErrorEmpiricoRecortado}. For convenience we state this
fact in the following result.
\begin{equation*}\label{DefEmpTrimErr}
R_{n,\alpha}(g):=\inf_{Q\in\mathcal{R}_\alpha (P_n)} Q(g(X)\ne Y),
\end{equation*}
where $P_n$ is the empirical distribution of $P$.

\begin{Corollary}\label{R_nalpha(g)}
Let $g$ be a given classifier, $\alpha$ a fixed trimming level
and $n\in\mathbb{N}$ the sample size,
\begin{equation}\label{CorrespErrorEmp}
R_{n,\alpha}(g)=\frac{1}{1-\alpha}(R_n(g)-\alpha)_+.
\end{equation}
\end{Corollary}

\medskip

In empirical risk minimization methods (see for instance
\cite{lugosi2002pattern} and references therein), the empirical
classification error $R_n(g)$ is used as an estimator of $R(g)$.
Among other good properties, $R_n(g)$ is unbiased as an estimator of
$R(g)$. This does not hold for trimmed errors. The following
proposition provides a control over this quantity and shows that the
empirical classification is still an asymptotically unbiased
estimate of the classification error.

\begin{proposition}\label{EsperanzaRnalpha}
For a given trimming level $\alpha$ and a given classifier $g$
\begin{equation*}
0\leq E(R_{n,\alpha}(g))- R_\alpha(g)\leq
\frac{\sqrt{R(g)}}{\sqrt{2n}(1-\alpha)}.
\end{equation*}
\end{proposition}

\section{Optimal selection of trimming levels in classification} \label{s:choice}
\subsection{Main Result}
Trimmed models enable to decrease the classification error in such a
way that the loss of information of using less observations can be
quantified and controlled. As in any robust procedure, we aim at
selecting the amount of data to be removed, which, in this setting,
corresponds to the optimal trimming level. Actually the aim is to
find a data driven $\hat{\alpha}$ such that the classification risk
is minimized without removing a too large quantity of information
about the initial distribution. We know that the bigger the trimming
is, the smaller the error will be, but the more data we trim, the
less information our model will keep. To look for an equilibrium we
will introduce a penalization which will depend on the size of the
chosen trimming level. For the sake of clarity we present first an
oracle bound in the toy setup in which we only consider a fixed
classification rule and we aim at choosing the right trimming
proportion. Later we present a more general result which will deal
with the more realistic case in which the classifier is chosen
within a more general collection of models.

\begin{theorem}\label{CorDesigOraculDisc}
Let $\xi_1=(Y_1,X_1),\ldots,\xi_n=(Y_n,X_n)$ be $n$  i.i.d
observations  with distribution P that take values in
$\{0,1\}\times\mathbb{R}^p$. Let $g$ be a given classifier and
$\alpha_{max}\in(0,1)$. Consider the penalization function
\begin{equation*}\label{FormulaPenCont}
pen(\alpha)= \frac{1}{(1-\alpha)}\sqrt{\frac{\ln(n)}{2n}}
\end{equation*}
and define
\begin{equation*}
\hat{\alpha}=\argmin_{\alpha\in [0,\alpha_{max}]}
R_{n,\alpha}(g)+pen(\alpha),
\end{equation*}
then the following bound holds,
\begin{equation*}\label{CotaPenaliaztionCont}
E(R_{\hat{\alpha}}(g))\leq \inf_{\alpha \in [0,\alpha_{max}]}
\left(R_\alpha(g)+pen(\alpha)+
\frac{\sqrt{R(g)}}{\sqrt{n}(1-\alpha)}\right)+\frac{1}{(1-\alpha_{max})}\sqrt{\frac{2\pi}{n}}+\frac{1}{n(1-\alpha_{max})^2}.
\end{equation*}
\end{theorem}

This theorem enables to understand the effect of trimming on the
classification error. For a given classifier $g$ we fix a maximum
level of trimming $\alpha_{max}$ that we do not want to exceed. Then
the automatic penalized rule for choosing the trimming level leads
to an oracle inequality that warrants that the best classification
error is achieved. Similar to model selection rules, the price to
pay is a term of order $1/\sqrt{n}$ which does not hamper the
classification error. In particular if the classifier $g$ has a
small classification error in the sense that $R(g)$ is smaller than
some $ \alpha <  \alpha_{\max}$, we achieve to remove the data that
are misclassified, leading to a smaller classification error. \vskip
.1in


A natural extension of this result is the case where we consider a
class of classification rules among which the optimal classifier
will be selected. A complex class will usually lead us to rules that
have a small bias in the sense that they classify well the data in
the training sample yet at the expense of  larger variance error,
usually leading to  an overfitting of the classification model. To
deal with this necessary control of complexity, the penalties will
not only depend on the trimming level as before but also on the
complexity of the class of classifiers. This complexity will be
measured using the Vapnik-Chervonenkis dimension (VC), see for
instant in~\cite{geer2000empirical} and references therein. Here
$\mathcal{F}$ denotes the set of all classifiers.

%

\begin{theorem} \label{CorDesigOraculDiscFuncFin}
Let $\xi_1,\ldots,\xi_n$ be $n$ independent and identically
distributed observations with distribution P that take values in
$\{0,1\}\times\mathbb{R}^p$. Let $\{\mathcal{G}_m\}_{m\in\mathbb{N}}
\subset \mathcal{F}$ be a family of classes of classifiers with
Vapnik-Chervonenkis dimension $V_{\mathcal{G}_m}<\infty$ for all
$m\in\mathbb{N}$. Let $\alpha_{max}\in(0,1)$ and let $\Sigma$ be a
non-negative constant. Consider $\{x_m\}_{m\in\mathbb{N}}$ a family
of non-negative weights such that
\begin{equation*}
\displaystyle\sum_{m\in\mathbb{N}}e^{-x_m}\leq\Sigma<\infty.
\end{equation*}
If we consider the penalization function
\begin{equation*}\label{FormulaPen}
pen(\alpha,\mathcal{G}_m)=\sqrt{\frac{\ln(n)+x_m}{2n(1-\alpha)^2}}+
\frac{1}{(1-\alpha)}\sqrt{\frac{V_{\mathcal{G}_m}\ln(n+1)+\ln(2)}{n}}
\end{equation*}
and we define
\begin{equation*}
(\hat{\alpha},\hat{m})=\argmin_{(\alpha,m)\in
[0,\alpha_{max}]\times\mathbb{N}}
R_{n,\alpha}(\mathcal{G}_m)+pen(\alpha,\mathcal{G}_m),
\end{equation*}
the following bound holds
\begin{align*}
E(R_{\hat{\alpha}}(\mathcal{G}_{\hat{m}})) & \leq \min_{(\alpha,m)
\in [0,\alpha_{max}]\times\mathbb{N}}
\left(R_\alpha(\mathcal{G}_m)+pen(\alpha,\mathcal{G}_m)+
  \frac{\sqrt{R(\mathcal{G}_m)}}{\sqrt{2n}(1-\alpha)}\right) \\ & +
\frac{1+\Sigma}{2(1-\alpha_{max})}\sqrt{\frac{\pi}{2n}}+
\frac{1}{n(1-\alpha_{max})^2}.
\end{align*}
\end{theorem}

Here, again, we obtain a bound similar to the result provided in
\cite{Massart2007libro}. The penalty for choosing the trimming
parameter depends on the VC dimension of the class of classifiers.
Hence, this choice leads to an oracle inequality ensuring the
optimality of this selection procedure. As before, the effect of
trimming is that it removes an optimal number of data that are
misclassified by the collection of classifiers, leading to better
classification rates on the set of \textit{good} data.

For a better understanding about the implications of Theorem
\ref{CorDesigOraculDiscFuncFin} we include next a section which
explores this bound for the particular case of linear classifiers.

\subsection{Example}

Assume we have $n$ i.i.d. observations $(Y_1,X_1),\ldots,(Y_n,X_n)$
where $X_i\in\mathbb{R}^p$ and $Y_i\in\{0,1\}$. We consider the
collection of models $\{\mathcal{G}_m\}_{m\in \mathcal{M}}$ where
for each $m$, $\mathcal{G}_m$ is the family of linear classifiers
built only using a selection of variables consisting of the first
$m$ components of $X_i$. Set $\mathcal{M}=\{1,\ldots,p\}$. For
$x\in\mathbb{R}^p$ let $x^{(m)}$ denote the vector consisting of the
first $m$ components of $x$. Define the set of possible classifiers
as
$$\mathcal{G}_m=\left\{g\in\mathcal{F}:g(x)=I_{[a^T x^{(m)}+b\geq0]};a\in\mathbb{R}^m;b\in\mathbb{R}\right\}.$$
Let us also denote by $\mathcal{A}_m$ the collection of all sets
\begin{equation*}
\{\{0\}\times \{x:g_m(x)=1\}\} \bigcup \{\{1\}\times
\{x:g_m(x)=0\}\}
\end{equation*}
and by $\mathcal{B}_m$ the collection of sets
$$\left\{x\in\mathbb{R}^p:g_m(x)=1\right\}$$
where $g_m$ ranges in $\mathcal{G}_m$. Using, for instance, Theorem
13.1 and Corollary 13.1 in \cite{LugosiPattern}, we have that
$V_{\mathcal{G}_m}=V_{\mathcal{A}_m}=m+1$. Then the penalization
function considered in Theorem \ref{CorDesigOraculDiscFuncFin} can
be written as
\begin{equation*}
pen(\alpha,\mathcal{G}_m)=\sqrt{\frac{\ln(n)+x_m}{2n(1-\alpha)^2}}+
\frac{1}{(1-\alpha)}\sqrt{\frac{(m+1)\ln(n+1)+\ln(2)}{n}}.
\end{equation*}

We will choose the family of non-negative weights $x_m=\ln(p)$ for
all $m\in\mathcal{M}$ and the universal constant $\Sigma=1$. If we
define
\begin{eqnarray*}
(\hat{\alpha},\hat{m})&=&\argmin_{(\alpha,m)\in
[0,\alpha_{max}]\times\mathcal{M}} \left(
R_{n,\alpha}(\mathcal{G}_m)+\sqrt{\frac{\ln(np)}{2n(1-\alpha)^2}} \right. \\
&+& \left. \frac{1}{(1-\alpha)}\sqrt{\frac{(m+1)\ln(n+1)+\ln(2)}{n}}
\right),
\end{eqnarray*}
this leads to the following bound
\begin{eqnarray*}
E(R_{\hat{\alpha}}(\mathcal{G}_{\hat{m}}))&\leq& \min_{(\alpha,m)
\in [0,\alpha_{max}]\times\mathcal{M}}
\left(R_\alpha(\mathcal{G}_m)+\sqrt{\frac{\ln(np)}{2n(1-\alpha)^2}}\right.\\
&+&\left. \frac{1}{(1-\alpha)}\sqrt{\frac{(m+1)\ln(n+1)+\ln(2)}{n}}+
\frac{\sqrt{R(\mathcal{G}_m)}}{\sqrt{2n}(1-\alpha)}\right)\\
&+& \frac{1}{(1-\alpha_{max})}\sqrt{\frac{\pi}{2n}}+
\frac{1}{n(1-\alpha_{max})^2}.
\end{eqnarray*}


First we point out that trimming reduces the classification error
that may vanish as long as $\alpha_{\max}$ is large enough to remove
a sufficient fraction of observations. As in model selection
techniques, the last three terms of the right hand side of the
inequality are of order $1/\sqrt{n}$ while for a fixed $m$ the third
term will be of order $\sqrt{m \ln(n)/n}$. Finally, the second term
is of order $\sqrt{\ln{(n)}/n+\ln{(p)}/n} $. Hence, as long as
$\ln(p)$ is smaller than $n$, the expected value of the best trimmed
classification error for the best class will be small as the number
of observations increases.

\section{Conclusions and perspectives} \label{s:conclu}

In classification theory, many classification rules are affected by
the presence of points which are very difficult to classify. When
dealing with high dimensional observations or when the number of
observations is large, this situation occurs quite often and may
drastically hamper the performance of classifiers which take into
account all the data. One may be tempted to focus on these points
and modify the classification rule to increase their classification
ranking for these special points. This is the point of view of
boosting algorithms for instance, as described in
\cite{freund1995desicion} for example. Yet this is often done at the
expense of the complexity of the rule and its ability to be
generalized. Hence a practical and maybe pragmatic solution is to
consider some of these points as outliers and to simply remove them.
Statisticians are reluctant to discard observations, yet in many
applications, in particular when confronted to large amounts of
observations, this enables to produce rules that are easier to
interpret and that can provide a better understanding of the
phenomenon which is studied, provided not too many data are removed
from the training sample. This is the typical choice made in several
papers but not much is said about the way the outliers are selected
and its impact on the classification performance.

This is the reason why we tried to provide in this paper a
statistical framework to robust classification by removal of some
observations. We provided a method that considers data as outliers
based on their classification error by a classifier or a given class
of classifiers. Within this framework this procedure enables to
select in a data driven way an optimal proportion of observations to
be removed in order to achieve a better classification error. The
level of trimming and the best classifier are selected
simultaneously and we obtain an oracle inequality to assess the
quality of this procedure. We think that this result may provide
some guidelines to remove outliers for classification problems with
theoretical guarantees.

Yet we rely on a minimization of a penalized $0-1$ loss function
which is difficult to handle. A version of this trimming procedure
for convex functions that lead to a feasible way of computing the
weights is actually under study. We will thus obtain a way of
detecting outliers and removing them such that the classification
error with this new data set will  be theoretically controlled.

\section{Appendix} \label{s:append}
\subsection{Technical lemmas}

Let $A\subset\{0,1\}\times\mathbb{R}^p$, we denote
$A_i=\{x\in\mathbb{R}^p:(i,x)\in A\}$, for $i=0,1$. Obviously
$A=(\{0\}\times A_0)\cup (\{1\}\times A_1)$ and the union is
disjoint, so for every measurable set
$A\subset\{0,1\}\times\mathbb{R}^p$ and every probability $P \in
\{0,1\}\times\mathbb{R}^p$,

\begin{equation}\label{notacionbasica}
    P(A)=p_0P_0(A_0)+p_1P_1(A_1),
\end{equation}
where $p_0=P(\{0\}\times\mathbb{R}^p)$, $p_1=1-p_0$,
$P_0(A_0)=P(A|Y=0)=P(\{0\}\times A_0)/p_0$ and
$P_1(A_1)=P(A|Y=1)=P(\{1\}\times A_1)/p_1$. $P_0$ and $P_1$ are
probabilities in $\mathbb{R}^p$. Conversely, from $p_0\in[0,1]$ and
the probabilities $P_0$ and $P_1$ in $\mathbb{R}^p$ the equation
\eqref{notacionbasica} defines a probability in
$\{0,1\}\times\mathbb{R}^p$ and the relation is one on one (except
for the degenerate cases in which $p_0=0$ or $p_0=1$), so we can
identify the probability P with the object $(p_0,P_0,P_1)$. We will
set $P \equiv (p_0,P_0,P_1)$.

\begin{Lemma}\label{recortesclas} With the previous notation, if $Q\equiv (q_0,Q_0,Q_1)$ with $q_0\in (0,1)$,
then $Q\in \mathcal{R}_\alpha (P)$ if and only if
\begin{equation}\label{equivrecortes}
q_0\leq \frac {p_0}{1-\alpha},\quad 1-q_0\leq \frac
{1-p_0}{1-\alpha}, \quad
Q_0\in\mathcal{R}_{1-\frac{q_0}{p_0}(1-\alpha)}(P_0)\quad \mbox{ and
}\quad Q_1\in\mathcal{R}_{1-\frac{1-q_0}{1-p_0}(1-\alpha)}(P_1).
\end{equation}
\end{Lemma}

\medskip
\noindent \textbf{Proof.} Note first that $q_0=Q(\{0\}\times
\mathbb{R}^p)$, $Q\in\mathcal{R}_\alpha (P)$ implies $q_0\leq
\frac{1}{1-\alpha} P(\{0\}\times  \mathbb{R}^p) =\frac
{p_0}{1-\alpha}$. The same argument shows that $1-q_0\leq \frac
{1-p_0}{1-\alpha}$ if $Q\in\mathcal{R}_\alpha (P)$. Observe that the
conditions $q_0\leq \frac {p_0}{1-\alpha}$ and $1-q_0\leq \frac
{1-p_0}{1-\alpha}$ guarantee that $0\leq
1-\frac{q_0}{p_0}(1-\alpha)\leq 1$ and $0\leq
1-\frac{1-q_0}{1-p_0}(1-\alpha)\leq 1$, hence the trimming sets in
the statement are well defined. Moreover, if $Q\in\mathcal{R}_\alpha
(P)$ then
$$Q_0(A_0)=\frac{Q(\{0\}\times A_0)}{q_0}\leq \frac{1}{(1-\alpha)q_0} P(\{0\}\times A_0)=\frac{1}{(1-\alpha)\frac{q_0}{p_0}} P_0(A_0), $$
which proves that
$Q_0\in\mathcal{R}_{1-\frac{q_0}{p_0}(1-\alpha)}(P_0)$. In a similar
way it can be proven that
$Q_1\in\mathcal{R}_{1-\frac{1-q_0}{1-p_0}(1-\alpha)}(P_1)$, which
proves that the assumptions~ \eqref{equivrecortes} are necessary. To prove  the sufficiency note that if we have \eqref{equivrecortes}
then $q_0Q_0(A_0)\leq \frac{1}{1-\alpha}P_0(A_0)$,
$(1-q_0)Q_1(A_1)\leq \frac{1}{1-\alpha}P_1(A_1)$ and hence
$$Q(A)=q_0Q_0(A_0)+(1-q_0)Q_1(A_1)\leq \frac 1 {1-\alpha}(p_0P_0(A_0)+(1-p_0)P_1(A_1))=\frac{1}{1-\alpha} P(A),$$
which completes the proof. \hfill $\Box$
\bigskip

With this identification we now prove the following lemma that will be the
first step to prove Proposition \ref{R_alpha(g)}.

\begin{Lemma}\label{LemaRRecT}
With the previous notation

\begin{equation}\label{RRecT}
R_\alpha(g)=\displaystyle\min_{1-\frac{1-p_0}{1-\alpha}\leq q_0\leq
\frac{p_0}{1-\alpha}} \left[
\left(q_0-\frac{p_0}{1-\alpha}P_0(g(x)=0)\right)_+ +
\left(1-q_0-\frac{1-p_0}{1-\alpha}P_1(g(x)=1)\right)_+ \right].
\end{equation}

\end{Lemma}

\medskip
\noindent \textbf{Proof.}

\noindent The first step consists in  writing the probability Q in
terms of  $(q_0,Q_0,Q_1)$

\begin{eqnarray}\label{ProbQ}
Q(g(x)\neq y) & = & q_0\displaystyle \int_{(g(x)=1)}
\frac{dQ_0}{d\mu} d\mu + (1-q_0)\displaystyle \int_{(g(x)=0)}
\frac{dQ_1}{d\mu} d\mu \nonumber\\
& = & \int \left(q_0 I_{(g(x)=1)}\frac{dQ_0}{d\mu} +
(1-q_0)I_{(g(x)=0)} \frac{dQ_1}{d\mu}\right) d\mu.
\end{eqnarray}

\noindent We are looking for the probability that minimizes the
probability of error between all the probabilities $Q\in
\mathcal{R}_\alpha (P)$, this means we are looking for $Q_0$ and
$Q_1$ that minimize \eqref{ProbQ}. We are going to make the
calculations for $Q_0$, $Q_1$ can be gotten analogously.

\noindent As we are minimizing, we are going to concentrate the
probability $Q_0$ in the set $(g(x)=0)$. By Lemma \ref{recortesclas}
we know that $Q_0\leq \frac{p_0}{q_0(1-\alpha)}P_0$, so the value of
$Q_0$ depends on the value of $P_0$. There are two possibilities,

\begin{enumerate}
  \item $P_0(g(x)=0)\geq \frac{q_0}{p_0}(1-\alpha)$: As
  $\frac{p_0}{q_0(1-\alpha)}P_0\geq1$ we can group all the probability $Q_0$ in the set $\{x\in\mathbb{R}^p /g(x)=0\}$ and hence $Q_0(g(x)=0)=1$.
  \item $P_0(g(x)=0) < \frac{q_0}{p_0}(1-\alpha)$: Now we can not give to $Q_0(g(x)=0)$ probability 1 because we would be violating the condition in Lemma \ref{recortesclas},
  hence $Q_0(g(x)=0)=\frac{P_0(g(x)=0)}{\frac{q_0}{p_0}(1-\alpha)}$.
\end{enumerate}

\noindent And in the optimum we will have

$$Q_0(g(x)=0)=\min\left(\frac{P_0(g(x)=0)}{\frac{q_0}{p_0}(1-\alpha)},1\right).$$

\noindent As we are concerned in $Q_0(g(x)=1)$ and $Q_0$ is a
distribution,
$$Q_0(g(x)=1)=\left(1-\frac{p_0}{q_0(1-\alpha)}P_0(g(x)=0)\right)_+,$$
analogously
$$Q_1(g(x)=0)=\left(1-\frac{1-p_0}{(1-q_0)(1-\alpha)}P_1(g(x)=1)\right)_+.$$

\noindent So for a fixed $q_0$ and $Q_0$, $Q_1$ as in Lemma
\ref{recortesclas}
\begin{eqnarray*}
  \displaystyle \min_{Q_0,Q_1} Q(g(x)\neq y) &=& q_0\left(1-\frac{p_0}{q_0(1-\alpha)}P_0(g(x)=0)\right)_+ \\
  &+& (1-q_0)\left(1-\frac{1-p_0}{(1-q_0)(1-\alpha)}P_1(g(x)=1)\right)_+.
\end{eqnarray*}

Using that $q_0$ and $1-q_0$ are positive lead to ~\eqref{RRecT}.
The limits for $q_0$ are obtained from Lemma \ref{recortesclas}.
\hfill $\Box$ \vskip .1in

We prove that both  the theoretical trimmed error  and the
empirical error of two close trimming levels is close.

\begin{proposition}\label{DifEntreErroresRecortados}
Let $\alpha_1$, $\alpha_2$ be two trimming levels such that
$\alpha_2 \in [\alpha_1,\alpha_1+\frac{1}{n}]$, let $\alpha_ {max}$
be such that $\alpha_1\leq\alpha_2\leq\alpha_{max}<1$ and let $g$ be
a given classifier, then
\begin{equation*}
R_{\alpha_1}(g)-R_{\alpha_2}(g)\leq \frac{1}{n(1-\alpha_{max})^2} \
\ \text{and}\ \ R_{n,\alpha_1}(g)-R_{n,\alpha_2}(g)\leq
\frac{1}{n(1-\alpha_{max})^2}.
\end{equation*}
\end{proposition}
\medskip
\noindent \textbf{Proof.}
\begin{align*}
&R_{\alpha_1}(g)-R_{\alpha_2}(g) =
\frac{(R(g)-\alpha_1)_+}{(1-\alpha_1)}-\frac{(R(g)-\alpha_2)_+}{(1-\alpha_2)}\\
&=\frac{((1-\alpha_2)(R(g)-\alpha_1))_+-((1-\alpha_1)(R(g)-\alpha_2))_+}{(1-\alpha_1)(1-\alpha_2)}\\
&\leq
\frac{1}{(1-\alpha_1)(1-\alpha_2)}|R(g)-\alpha_1-\alpha_2R(g)+\alpha_1\alpha_2-(R(g)-\alpha_2-\alpha_1R(g)+\alpha_1\alpha_2)|\\
&=\frac{1}{(1-\alpha_1)(1-\alpha_2)}|-\alpha_1-\alpha_2R(g)+\alpha_2+\alpha_1R(g)|\\
&=\frac{1}{(1-\alpha_1)(1-\alpha_2)}|(R(g)-1)(\alpha_1-\alpha_2)|=\frac{1}{(1-\alpha_1)(1-\alpha_2)}|R(g)-1||\alpha_1-\alpha_2|.
\end{align*}
As we chose $\alpha_2$, $|\alpha_1-\alpha_2|\leq\frac{1}{n}$ and as
for every value of $\alpha$ we can bound $\frac{1}{1-\alpha}$ by
$\frac{1}{1-\alpha_{max}}$ and $|R(g)-1|\leq 1$, we can conclude
that
\begin{equation*}
R_{n,\alpha_1}(g)-R_{n,\alpha_2}(g)\leq
\frac{1}{n(1-\alpha_{max})^2}.
\end{equation*}
\noindent The proof is identical for the empirical trimmed error.
\hfill $\Box$

\subsection{Proofs}
\noindent \textbf{Proof of Proposition \ref{R_alpha(g)}.} The result
is a direct consequence of the  minimization with respect to  $q_0$
of the expression obtained in Lemma~\ref{LemaRRecT}.

First  see that $R_\alpha(g)=0$ if and
only if $1-\frac{1-p_0}{1-\alpha}P_1(g(x)=1)\leq
\frac{p_0}{1-\alpha}P_0(g(x)=0)$. Then consider
the opposite case.

\noindent As we are adding two positive terms, the sum is equal to 0
only if  both terms are equal to 0, leading to
$$\left(q_0-\frac{p_0}{1-\alpha}P_0(g(x)=0)\right)_+\leq 0 \
\Leftrightarrow \ q_0 \leq \frac{p_0}{1-\alpha}P_0(g(x)=0),$$ in a
similar way we obtain
$$\left(1-q_0-\frac{1-p_0}{1-\alpha}P_1(g(x)=1)\right)_+\leq 0 \
\Leftrightarrow \ q_0 \geq 1-\frac{1-p_0}{1-\alpha}P_1(g(x)=1).$$ So
$R_\alpha(g)=0$ if and only if
$1-\frac{1-p_0}{1-\alpha}P_1(g(x)=1)\leq
\frac{p_0}{1-\alpha}P_0(g(x)=0)$.

\noindent Now consider the case where this inequality does not hold,
this means, $1-\frac{1-p_0}{1-\alpha}P_1(g(x)=1)>
\frac{p_0}{1-\alpha}P_0(g(x)=0)$. The first term of \eqref{RRecT} is
a stepwise lineal function with value $0$ until
$\frac{p_0}{1-\alpha}P_0(g(x)=0)$ and increasing with slope $1$
since then. The second term is also stepwise linear, in this case it
decreases with slope $-1$ until it reaches $0$ in
$1-\frac{1-p_0}{1-\alpha}P_1(g(x)=1)$ with  value  $0$ from that
point.

Now we are going to see that in this case the interval
$\left[\frac{p_0}{1-\alpha}P_0(g(x)=0),1-\frac{1-p_0}{1-\alpha}P_1(g(x)=1)\right]$
gives us the minimal value of \eqref{RRecT}. If
$1-\frac{1-p_0}{1-\alpha}P_1(g(x)=1)<1-\frac{1-p_0}{1-\alpha}$ or
$\frac{p_0}{1-\alpha}P_0(g(x)=0)>\frac{p_0}{1-\alpha}$ we will
eliminate non feasible values of $q_0$ from the optimal set, hence
the set of optimal values of $q_0$ that minimizes $R_\alpha(g)$ are
\begin{equation}\label{q0IntOptim}
q_0(x)=\left\{ \begin{array}{lcc}
             \left[1-\frac{1-p_0}{1-\alpha}P_1(g(x)=1),\frac{p_0}{1-\alpha}P_0(g(x)=0)\right] \bigcap\left[1-\frac{1-p_0}{1-\alpha},\frac{p_0}{1-\alpha}\right] &   si  &R(g)\leq \alpha \\
             \\ \left[\frac{p_0}{1-\alpha}P_0(g(x)=0),1-\frac{1-p_0}{1-\alpha}P_1(g(x)=1)\right] \bigcap \left[1-\frac{1-p_0}{1-\alpha},\frac{p_0}{1-\alpha}\right] &  si & R(g)>\alpha \\
             \end{array}
   \right..
\end{equation}

\noindent Let us see this. We are going to suppose, for simplicity,
that we are in the case
$$1-\frac{1-p_0}{1-\alpha} \leq \frac{p_0}{1-\alpha}P_0(g(x)=0) <
1-\frac{1-p_0}{1-\alpha}P_1(g(x)=1) \leq \frac{p_0}{1-\alpha}.$$ Let
$I_1=[1-\frac{1-p_0}{1-\alpha},\frac{p_0}{1-\alpha}P_0(g(x)=0)]$,
$I_2=[\frac{p_0}{1-\alpha}P_0(g(x)=0),1-\frac{1-p_0}{1-\alpha}P_1(g(x)=1)]$,
$I_3=[1-\frac{1-p_0}{1-\alpha}P_1(g(x)=1),\frac{p_0}{1-\alpha}]$, we
denote
$$R^i=\displaystyle\min_{I_i} \left[
\left(q_0-\frac{p_0}{1-\alpha}P_0(g(x)=0)\right)_+ +
\left(1-q_0-\frac{1-p_0}{1-\alpha}P_1(g(x)=1)\right)_+ \right],$$
for $i=1,2,3$. Obviously $R_\alpha(g)=\min R^i$.

\noindent In $I_1$ the first term is 0 because $q_0\leq
\frac{p_0}{1-\alpha}P_0(g(x)=0)$ and the second term is
$1-\frac{1-p_0}{1-\alpha}P_1(g(x)=1)-q_0$. As we are looking for a
minimization of this value and $q_0$ is subtracting, we will give to
it the biggest value it can take, that is, the upper bound of the
interval. Hence,
\begin{eqnarray*}
R^1 & = &
1-\frac{1-p_0}{1-\alpha}P_1(g(x)=1)-\frac{p_0}{1-\alpha}P_0(g(x)=0)\\
& = &
1-\frac{(1-p_0)(1-P_1(g(x)=0))+p_0(1-P_0(g(x)=1))}{1-\alpha}\\
& = & 1-\frac{1-R(g)}{1-\alpha}.
\end{eqnarray*}

\noindent If we are in $I_2$ none of the terms is going to be 0.
First one is $q_0-\frac{p_0}{1-\alpha}P_0(g(x)=0)$ and second one
$1-\frac{1-p_0}{1-\alpha}P_1(g(x)=1)-q_0$, when we add them, the
$q_0$ in both terms clears and we obtain
$$R^2=1-\frac{p_0}{1-\alpha}P_0(g(x)=0)-\frac{1-p_0}{1-\alpha}P_1(g(x)=1)=1-\frac{1-R(g)}{1-\alpha}.$$

\noindent Last, in $I_3$ is the second term which becomes 0, letting
the first one as $q_0-\frac{p_0}{1-\alpha}P_0(g(x)=0)$. In this case
$q_0$ is adding, we want to give the minimum value possible so
$$R^3=1-\frac{1-p_0}{1-\alpha}P_1(g(x)=1)-\frac{p_0}{1-\alpha}P_0(g(x)=0)=1-\frac{1-R(g)}{1-\alpha}.$$

\noindent And, as we have already said, the minimum is attained at
$$\left[1-\frac{1-p_0}{1-\alpha}P_1(g(x)=1),
\frac{p_0}{1-\alpha}P_0(g(x)=0)\right].$$ Moreover, since
$R^1=R^2=R^3$, the value of this minimum will be
$$R_\alpha(g)=1-\frac{1-R(g)}{1-\alpha}.$$

\noindent Putting together both cases we have that $R_\alpha(g)$
reaches its minimum in \eqref{q0IntOptim} and, since condition
$1-\frac{1-p_0}{1-\alpha}P_1(g(x)=1)>\frac{p_0}{1-\alpha}P_0(g(x)=0)$
holds if and only if $R(g)>\alpha$, we have that $R_\alpha(g)=0 \
\Leftrightarrow \ R(g)\leq \alpha$ and hence,
$$R_\alpha(g)=\frac{1}{1-\alpha}(R(g)-\alpha)_+.$$ \hfill $\Box$

\medskip
\noindent \textbf{Proof of Proposition~\ref{RelBayesrec}}. Note that
$Err(P)=R(g_B)$ and $Err_\alpha(P)=R_\alpha(g_B^\alpha)$. Recall
that
\begin{eqnarray*}
\mbox{Err}_\alpha(P)&:=& \inf_{Q\in\mathcal{R}_\alpha
(P)}\mbox{Err}(Q)=\inf_{Q\in\mathcal{R}_\alpha (P)} \inf_g
Q(g(x)\neq y)= \inf_g \inf_{Q\in\mathcal{R}_\alpha (P)} Q(g(x)\neq
y)\\
&=&\inf_g R_\alpha(g)=\min_g\frac{(R(g)-\alpha)_+}{1-\alpha},
\end{eqnarray*}
the minimum in the last inequality is due to Proposition
\ref{R_alpha(g)}. The infimum is reached so it is a minimum.
Moreover we know that this error is minimal when the classifier is
Bayes classifier, so
$$\mbox{Err}_\alpha(P)=\frac{(R(g_B)-\alpha)_+}{1-\alpha}=\frac{(\mbox{Err}(P)-\alpha)_+}{1-\alpha}.$$
\hfill $\Box$

\medskip
\noindent \textbf{Proof of Proposition~\ref{EsperanzaRnalpha}.} The
first inequality can be proved by
\begin{eqnarray*}
E(R_{n,\alpha}(g)) &=&
E\left(\frac{(R_n(g)-\alpha)_+}{1-\alpha}\right)=\frac{1}{1-\alpha}E\left((R_n(g)-\alpha)_+\right)\\
&\geq&
\frac{1}{1-\alpha}\left(E(R_n(g))-\alpha\right)_+=\frac{1}{1-\alpha}\left(R(g)-\alpha\right)_+=R_\alpha(g).
\end{eqnarray*}
\noindent Where we have used \eqref{CorrespErrorEmp} for the first
equality and the property $E(R_n(g))=R(g)$ and \eqref{RalphaFunR}
for the two last ones. The inequality comes from applying Jensen
inequality, and this is possible due to the fact that $(.)_+$ is a
convex function.

\noindent For the second inequality we need Proposition
\ref{R_alpha(g)} and by Corollary \ref{R_nalpha(g)},
\begin{equation}\label{LemaRelERnyRIneq1}
E(R_{n,\alpha}(g))-R_\alpha(g) =
\frac{E((R_n(g)-\alpha)_+)-(R(g)-\alpha)_+}{(1-\alpha)}.
\end{equation}
Let $X$ be a random variable such that $X=R_n(g)$, we know from
\cite{lugosi2002pattern} that $E(X)=R(g)$, and let
$\varphi(x)=(x-\alpha)_+$. $\varphi$ is a convex function, so
Jensen's inequality can be applied, this means $\varphi(E(X))\leq
E(\varphi(X))$. This function also is 1-Lispchitz and increasing, so
it satisfies the property $\varphi(y)-\varphi(x)\leq(y-x)_+$.

\noindent As we are not modifying $\frac{1}{1-\alpha}$ we are going
to let it aside for the moment and we will focus in the numerator,
applying $X$'s definition, mean's properties and $\varphi$'s
property mentioned above
\begin{eqnarray*}
E((R_n(g)-\alpha)_+)-(R(g)-\alpha)_+&=&E(\varphi(X))-\varphi(E(X))=E(\varphi(X)-\varphi(E(X)))\\
&\leq& E((X-E(X))_+).
\end{eqnarray*}
Now let $Y$ be a random variable such that $Y=^d X$, $Y$ and $X$ are
independent, this implies $E(Y)=E_X(Y)$, applying this, again mean's
properties, Jensen's inequality (for $(.)_+$) and conditional mean's
properties we get
\begin{eqnarray*}
E((X-E(X))_+)&=&E((X-E(Y))_+)=E((X-E_X(Y))_+)=E((E_X(X-Y))_+)\\
&\leq& E(E_X((X-Y)_+))=E((X-Y)_+).
\end{eqnarray*}
Now we are using that $X-Y$ is a symmetric variable, that it also is
a centered variable, variance's property for the sum of two
independent variables and that $X$ and $Y$ are identically
distributed to obtain
\begin{eqnarray*}
E((X-Y)_+)&=&\frac{1}{2}E(X-Y)\leq\frac{1}{2}(Var(X-Y))^{1/2}=\frac{1}{2}(Var(X)+Var(Y))^{1/2}\\
&=&\frac{1}{2}(2Var(X))^{1/2}.
\end{eqnarray*}
Last we are using variance's properties, the fact that $nX\sim
b(n,R(g))$ and that $(1-R(g))\leq1$ and we obtain
\begin{equation*}
\frac{1}{2}(2Var(X))^{1/2}=\frac{1}{2}(2\frac{1}{n^2}Var(nX))^{1/2}=\frac{1}{2}(2\frac{1}{n^2}nR(g)(1-R(g)))^{1/2}=\frac{1}{\sqrt{2n}}\sqrt{R(g)}.
\end{equation*}
Joining this with \eqref{LemaRelERnyRIneq1}, we get
\begin{equation*}
E(R_{n,\alpha}(g))-R_\alpha(g)\leq
\frac{\sqrt{R(g)}}{\sqrt{2n}(1-\alpha)}.
\end{equation*}
\hfill $\Box$ \vskip .1in

\medskip
\noindent \textbf{Proof of Theorem \ref{CorDesigOraculDisc}.}
\noindent First consider the case where the trimming parameter only takes value in the discrete set $A=[0,\frac{1}{n},\ldots,\frac{k_0}{n}]$ with
$k_0=[n\alpha_{max}]$. By definition $\hat{\alpha}$ satisfies that
for all $ \alpha \in A$
\begin{equation*}\label{EqBasicaTeoCota1}
R_{n,\hat{\alpha}}(g)+pen(\hat{\alpha})\leq
R_{n,\alpha}(g)+pen(\alpha).
\end{equation*}
\noindent This implies that
\begin{equation*}
R_{\hat{\alpha}}(g)-R_{\hat{\alpha}}(g)+R_{n,\hat{\alpha}}(g)+pen(\hat{\alpha})\leq
R_{\alpha}(g)-R_{\alpha}(g)+R_{n,\alpha}(g)+pen(\alpha)
\end{equation*}
or, what is the same
\begin{equation}\label{EqBasicaTeoCota2}
R_{\hat{\alpha}}(g) \leq R_{\alpha}(g) + pen(\alpha) +
(R_{n,\alpha}(g)-R_{\alpha}(g)) - pen(\hat{\alpha}) +
(R_{\hat{\alpha}}(g)-R_{n,\hat{\alpha}}(g)).
\end{equation}

\noindent Let us focus in the inside of the parenthesis,
\begin{equation*}
R_{n,\alpha}(g)-R_{\alpha}(g)=[R_{n,\alpha}(g)-E(R_{n,\alpha}(g))]+[E(R_{n,\alpha}(g))-R_{\alpha}(g)],
\end{equation*}
by Proposition~\ref{EsperanzaRnalpha} the second brace can be
bounded by $\frac{\sqrt{R(g)}}{\sqrt{2n}(1-\alpha)}$. For first
brace we will apply McDiarmid's inequality  taking
$R_{n,\alpha}(g)=F(\xi_1,\ldots,\xi_n)$ where $\xi_i=(Y_i,X_i)$. As
\begin{equation*}
|F(\xi_1,\ldots,\xi_i,\ldots,\xi_n)-F(\xi_1,\ldots,\xi'_i,\ldots,\xi_n)|\leq\frac{1}{n(1-\alpha)},
\end{equation*}
we can apply the inequality and hence
\begin{equation*}\label{McDiarmid}
P(R_{n,\alpha}(g)-E(R_{n,\alpha}(g))\geq t)\leq
e^{-2t^2n(1-\alpha)^2}.
\end{equation*}
Given $z>0$ take $t=\sqrt{\frac{z}{2n(1-\alpha)^2}}$, we get
\begin{equation*}
P\left(R_{n,\alpha}(g)-E(R_{n,\alpha}(g))\geq
\sqrt{\frac{z}{2n(1-\alpha)^2}}\right)\leq e^{-z}.
\end{equation*}
Joining this with \eqref{EqBasicaTeoCota2}, we get that, except in a
set of probability not greater than $e^{-z}$,
\begin{equation}\label{EqBasicaTeoCota3}
R_{\hat{\alpha}}(g) \leq R_{\alpha}(g) + pen(\alpha) +
\frac{\sqrt{R(g)}}{\sqrt{2n}(1-\alpha)}
+\sqrt{\frac{z}{2n(1-\alpha)^2}} - pen(\hat{\alpha}) +
(R_{\hat{\alpha}}(g)-R_{n,\hat{\alpha}}(g)).
\end{equation}

\noindent We are going to focus now in the other parenthesis. As we
saw in Proposition \ref{EsperanzaRnalpha}
\begin{equation*}
R_{\hat{\alpha}}(g)-R_{n,\hat{\alpha}}(g)\leq \displaystyle
\sup_{\alpha\in A} (R_{\alpha}(g)-R_{n,\alpha}(g))\leq \displaystyle
\sup_{\alpha\in A} (E(R_{n,\alpha}(g))-R_{n,\alpha}(g)).
\end{equation*}
Applying again McDiarmid's inequality taking this time
$t=\sqrt{\frac{\ln(n)+z}{2n(1-\alpha)^2}}$ we have $\forall \alpha'
\in A$
\begin{equation*}\label{McDiarmidTeoDisc}
P\left(E(R_{n,\alpha'}(g))-R_{n,\alpha'}(g) \geq
\sqrt{\frac{\ln(n)+z}{2n(1-\alpha')^2}}\right)\leq
\frac{1}{n}e^{-z}.
\end{equation*}
As we were interested in calculating this probability for
$\hat{\alpha}$ we have
\begin{align*}
&P\left(\displaystyle \sup_{\alpha\in A}
(E(R_{n,\alpha}(g))-R_{n,\alpha}(g)) \geq
\sqrt{\frac{\ln(n)+z}{2n(1-\hat{\alpha})^2}}\right)\\
&\leq \displaystyle\sum_{\alpha'\in A}
P\left(E(R_{n,\alpha'}(g))-R_{n,\alpha'}(g) \geq
\sqrt{\frac{\ln(n)+z}{2n(1-\alpha')^2}}\right)\\
&\leq n\frac{1}{n}e^{-z}\leq e^{-z}.
\end{align*}

\noindent Hence with probability at least $1-e^{-z}$
\begin{equation} \label{CondEspTeoCotaCont}
E(R_{n,\hat{\alpha}}(g))-R_{n,\hat{\alpha}}(g) \leq
\sqrt{\frac{\ln(n)+z}{2n(1-\hat{\alpha})^2}}.
\end{equation}

\noindent Now let us consider the complete interval. If
$\alpha'\in[0,\alpha_{max}]$ there exists $\alpha''\in A$ such that
$\alpha''\leq\alpha'\leq\alpha''+\frac{1}{n}$. Then by Proposition
\ref{DifEntreErroresRecortados}, in the set where
\eqref{CondEspTeoCotaCont} is satisfied we have

\begin{align*}
&E(R_{n,\alpha'}(g))-R_{n,\alpha'}(g)\\
&=E(R_{n,\alpha''}(g))-R_{n,\alpha''}(g)+E(R_{n,\alpha'}(g)-R_{n,\alpha''}(g))-(R_{n,\alpha'}(g)-R_{n,\alpha''}(g))\\
&\leq
\sqrt{\frac{\ln(n)+z}{2n(1-\alpha'')^2}}+\frac{1}{n(1-\alpha_{max})^2}
\leq
\sqrt{\frac{\ln(n)+z}{2n(1-\alpha')^2}}+\frac{1}{n(1-\alpha_{max})^2}\\
&\leq\sqrt{\frac{\ln(n)}{2n(1-\alpha')^2}}+\sqrt{\frac{z}{2n(1-\alpha')^2}}+\frac{1}{n(1-\alpha_{max})^2}
\end{align*}
for all $\alpha'\in[0,\alpha_{max}]$.

\noindent If we take as penalty
\begin{equation*}
pen(\alpha) = \sqrt{\frac{\ln(n)}{2n(1-\alpha)^2}}
\end{equation*}
and we substitute in \eqref{EqBasicaTeoCota3}. Given that, by
Proposition \ref{EsperanzaRnalpha}, $R_{\alpha'}(g)\leq
E(R_{n,\alpha'}(g))$, with probability at least $1-2e^{-z}$
\begin{align*}
&R_{\hat{\alpha}}(g)\\
&\leq R_{\alpha}(g) + pen(\alpha) +
\frac{\sqrt{R(g)}}{\sqrt{2n}(1-\alpha)}
+\sqrt{\frac{z}{2n(1-\alpha)^2}} - pen(\hat{\alpha}) +
\sqrt{\frac{\ln(n)+z}{2n(1-\alpha')^2}}+\frac{1}{n(1-\alpha_{max})^2}\\
&\leq R_{\alpha}(g) + pen(\alpha) +
\frac{\sqrt{R(g)}}{\sqrt{2n}(1-\alpha)}
+2\sqrt{\frac{z}{2n(1-\alpha_{max})^2}}+\frac{1}{n(1-\alpha_{max})^2}.
\end{align*}
\noindent Integrating with respect to z and taking the infimum for
$\alpha\in [0,\alpha_{max}]$ we can conclude that
\begin{equation*}
E(R_{\hat{\alpha}}(g))\leq \inf_{\alpha \in [0,\alpha_{max}]}
\left(R_\alpha(g)+pen(\alpha)+
\frac{\sqrt{R(g)}}{\sqrt{n}(1-\alpha)}\right)+\frac{1}{(1-\alpha_{max})}\sqrt{\frac{2\pi}{n}}+\frac{1}{n(1-\alpha_{max})^2}.
\end{equation*}
\hfill $\Box$ \vskip .1in

\medskip
\noindent \textbf{Proof of Theorem \ref{CorDesigOraculDiscFuncFin}.}
To proof the theorem we need the following elemental result whose
proof is omitted.

\begin{Lemma}\label{ResultElemsupymin} Given two functions $f$ and $g$
and a real number $k>0$,
\begin{equation*}
|f(x)-g(x)|\leq k \ \Rightarrow \ |\displaystyle\sup_x
f(x)-\displaystyle\sup_x g(x)|\leq k,
\end{equation*}

\begin{equation*}
|f(x)-g(x)|\leq k \ \Rightarrow \ |\displaystyle\min_x
f(x)-\displaystyle\min_x g(x)|\leq k.
\end{equation*}
\end{Lemma}

\noindent As in the previous Theorem we define the set
$A=\{0,\frac{1}{n},\frac{2}{n},\ldots\frac{k_0}{n}\}$ with
$k_0=[n\alpha_{max}]$. Then, by definition, $\hat{\alpha}$ and
$\hat{m}$ satisfy that for all $ \alpha \in A$ and $ m\in\mathbb{N}$
\begin{equation*}\label{EqBasicaTeoCota1FuncFin}
R_{n,\hat{\alpha}}(\mathcal{G}_{\hat{m}})+pen(\hat{\alpha},\mathcal{G}_{\hat{m}})\leq
R_{n,\alpha}(\mathcal{G}_m)+pen(\alpha,\mathcal{G}_m).
\end{equation*}
Adding and subtracting $R_{\hat{\alpha}}(\mathcal{G}_m)$ and
$R_\alpha(\mathcal{G}_m)$ and organizing the terms we get the
following inequality. We  bound the remaining terms in the
parenthesis,

\begin{equation*}
R_{\hat{\alpha}}(\mathcal{G}_{\hat{m}}) \leq
R_{\alpha}(\mathcal{G}_m) + pen(\alpha,\mathcal{G}_m) +
(R_{n,\alpha}(\mathcal{G}_m)-R_{\alpha}(\mathcal{G}_m)) -
pen(\hat{\alpha},\mathcal{G}_{\hat{m}}) +
(R_{\hat{\alpha}}(\mathcal{G}_{\hat{m}})-R_{n,\hat{\alpha}}(\mathcal{G}_{\hat{m}})).
\end{equation*}

\noindent First we are going to bound
\begin{equation*}
R_{n,\alpha}(\mathcal{G}_m)-R_{\alpha}(\mathcal{G}_m)=
\displaystyle\min_{g\in\mathcal{G}_m} R_{n,\alpha}(\mathcal{G}_m)-
\displaystyle\min_{g\in\mathcal{G}_m} R_{\alpha}(\mathcal{G}_m)\leq
R_{n,\alpha}(g')-R_{\alpha}(g')
\end{equation*}
with $g':=\displaystyle\argmin_{g\in\mathcal{G}_m}R_\alpha(g)$. We
are now in the same conditions as in Theorem
\ref{CorDesigOraculDisc} and we can bound these quantities except on
a set of probability not greater than $e^{-z}$, with a given $z>0$
by
\begin{equation*}
R_{n,\alpha}(g')-R_{\alpha}(g')\leq\frac{R(g')}{\sqrt{2n}(1-\alpha)}+\sqrt{\frac{z}{2n(1-\alpha)^2}},
\end{equation*}
which leads us to
\begin{equation}\label{EqBasicaTeoCota3FuncFin}
R_{\hat{\alpha}}(\mathcal{G}_{\hat{m}}) \leq
R_{\alpha}(\mathcal{G}_m) + pen(\alpha,\mathcal{G}_m) +
\frac{\sqrt{R(g')}}{\sqrt{2n}(1-\alpha)}
+\sqrt{\frac{z}{2n(1-\alpha)^2}} -
pen(\hat{\alpha},\mathcal{G}_{\hat{m}}) +
(R_{\hat{\alpha}}(\mathcal{G}_{\hat{m}})-R_{n,\hat{\alpha}}(\mathcal{G}_{\hat{m}})).
\end{equation}

\noindent Now we want to bound
\begin{equation*}
R_{\hat{\alpha}}(\mathcal{G}_{\hat{m}})-R_{n,\hat{\alpha}}(\mathcal{G}_{\hat{m}})
\leq \displaystyle\sup_{(\alpha',m')\in A\times\mathbb{N}}
(R_{\alpha'}(\mathcal{G}_{m'})-R_{n,\alpha'}(\mathcal{G}_{m'})) \leq
\displaystyle\sup_{(\alpha',m')\in
A\times\mathbb{N}}\displaystyle\sup_{g\in\mathcal{G}_{m'}}
(R_{\alpha'}(g)-R_{n,\alpha'}(g)).
\end{equation*}

\noindent Let us focus on

\begin{eqnarray}
\label{TeoDesOrA1}\displaystyle\sup_{g\in\mathcal{G}_m}
(R_{\alpha}(g)-R_{n,\alpha}(g))&=&
E\left(\displaystyle\sup_{g\in\mathcal{G}_m}
(R_{\alpha}(g)-R_{n,\alpha}(g))\right)\\
\label{TeoDesOrA2}&+&\left[\displaystyle\sup_{g\in\mathcal{G}_m}
(R_{\alpha}(g)-R_{n,\alpha}(g))
-E\left(\displaystyle\sup_{g\in\mathcal{G}_m}
(R_{\alpha}(g)-R_{n,\alpha}(g))\right)\right].
\end{eqnarray}

\noindent To bounding \eqref{TeoDesOrA2} we will make use again of McDiarmid's
inequality. First we need to see that the bounded difference
conditions is met.

\noindent We define
$Z:=f(\xi_1,\ldots,\xi_n)=\displaystyle\sup_{g\in\mathcal{G}_m}
(R_{\alpha}(g)-R_{n,\alpha}(g))$ and
$Z^{(i)}:=f(\xi_1,\ldots,\xi'_i,\ldots,\xi_n)$, we want to prove
\begin{equation}\label{ConDifAcotTeoFunc}
|Z-Z^{(i)}|\leq c_i,
\end{equation}
for certain constants $c_i$. The empirical error $R_{n,\alpha}(g)$
is defined as in \eqref{DefErrorEmpiricoRecortado} and
$R^{(i)}_{n,\alpha}(g)$ as the empirical error associated to the
sample $\xi_1,\ldots,\xi'_i,\ldots,\xi_n$. We start from
\begin{equation*}
|(R_{\alpha}(g)-R_{n,\alpha}(g))-
(R_{\alpha}(g)-R^{(i)}_{n,\alpha}(g))|
\end{equation*}
which implies, using Lemma \ref{ResultElemsupymin}, that
\eqref{ConDifAcotTeoFunc}.
\begin{equation*}
|R_{n,\alpha}(g)-R^{(i)}_{n,\alpha}(g)|=\left|\displaystyle\min_{(w_1,\ldots,w_n)}\sum_j
w_jI_{(g(X_j)\neq Y_j)}-\displaystyle\min_{(w_1,\ldots,w_n)}\sum_j
w_jI_{(g(X'_j)\neq Y'_j)}\right|,
\end{equation*}
where $(Y',X')$ stands for the sample
$\xi_1,\ldots,\xi'_i,\ldots,\xi_n$. For a vector $(w_1,\ldots,w_n)$
that satisfies the conditions \eqref{WeightsConditions},
\begin{equation*}
\left|\displaystyle\sum_j w_jI_{(g(X_j)\neq
Y_j)}-\displaystyle\sum_j w_jI_{(g(X'_j)\neq
Y'_j)}\right|=w_j\left|(I_{(g(X_i)\neq Y_i)}-I_{g(X'_i)\neq
Y'_i})\right|\leq\frac{1}{n(1-\alpha)}.
\end{equation*}
And using the second statement of Lemma \ref{ResultElemsupymin} leads to
\begin{equation*}
|R_{n,\alpha}(g)-R^{(i)}_{n,\alpha}(g)|\leq\frac{1}{n(1-\alpha)},
\end{equation*}
or written in a different way
\begin{equation*}
|(R_{\alpha}(g)-R_{n,\alpha}(g))-
(R_{\alpha}(g)-R^{(i)}_{n,\alpha}(g))|\leq\frac{1}{n(1-\alpha)}.
\end{equation*}
Applying again Lemma \ref{ResultElemsupymin}, we get to
\eqref{ConDifAcotTeoFunc} with $c_i=\frac{1}{n(1-\alpha)}$. Now we
can use McDiarmid's inequality to prove
\begin{equation}\label{McDiarSupR-Rn}
P\left(\displaystyle\sup_{g\in\mathcal{G}_m}
(R_{\alpha}(g)-R_{n,\alpha}(g))
-E(\displaystyle\sup_{g\in\mathcal{G}_m}
(R_{\alpha}(g)-R_{n,\alpha}(g)))\geq\sqrt{\frac{\ln(n)+z+x_m}{2n(1-\alpha)^2}}\right)\leq\frac{1}{n}e^{-z-x_m}.
\end{equation}

\noindent To bound \eqref{TeoDesOrA1} we will use
Vapnik-Chervonenkis theory from \cite{LugosiCombMeth} or
\cite{lugosi2002pattern}. Before we are able to apply this theory we
need to transform our functions in suitable functions. For this we
will use the equalities \eqref{RalphaFunR} and
\eqref{CorrespErrorEmp} and the fact that the function positive
part, defines as $X_+:=\max(0,X)$, is Lipschitz.
\begin{eqnarray}\label{CotaEspDesigVapnik}
\nonumber E\left(\sup_{g\in\mathcal{G}_m}(R_\alpha(g)-R_{n,\alpha}(g))\right)&=&\frac{1}{1-\alpha}E\left(\sup_{g\in\mathcal{G}_m}((R(g)-\alpha)_+(R_n(g)-\alpha)_+)\right)\\
\nonumber &\leq&\frac{1}{1-\alpha}E\left(\sup_{g\in\mathcal{G}_m}|R(g)-R_n(g)|\right)\\
&\leq&\frac{2}{1-\alpha}\sqrt{\frac{V_{\mathcal{A}_{\mathcal{G}_m}}\ln(n+1)+\ln(2)}{n}}.
\end{eqnarray}
The last inequality comes from section 4.2 in \cite{LugosiCombMeth}.
Joining \eqref{McDiarSupR-Rn} and \eqref{CotaEspDesigVapnik} we get
$\forall \alpha'\in A$ and $\forall m' \in \mathbb{N}$
\begin{equation}\label{EqIniCorolarioCont}
P\left(\sup_{g\in\mathcal{G}_{m'}}(R_{\alpha'}(g)-R_{n,\alpha'}(g))
\geq \sqrt{\frac{\ln(n)+z+x_{m'}}{2n(1-\alpha')^2}}
+\frac{2}{1-\alpha'}\sqrt{\frac{V_{\mathcal{G}_{m'}}\ln(n+1)+\ln(2)}{n}}\right)\leq\frac{1}{n}e^{-z-x_{m'}}.
\end{equation}

\noindent As we are looking for a bound for
$R_{\hat{\alpha}}(\mathcal{G}_{\hat{m}})-R_{n,\hat{\alpha}}(\mathcal{G}_{\hat{m}})$,
we have that
\begin{align*}
&P\left(\bigcup_{(\alpha',m')\in{A\times\mathbb{N}}}\sup_{g\in\mathcal{G}_{m'}}(R_{\alpha'}(g)-R_{n,\alpha'}(g))
\geq \sqrt{\frac{\ln(n)+z+x_{m'}}{2n(1-\alpha')^2}}
+\frac{2}{1-\alpha'}\sqrt{\frac{V_{\mathcal{G}_{m'}}\ln(n+1)+\ln(2)}{n}}\right)\\
&\leq \displaystyle\sum_{\alpha'\in A}\sum_{m'\in\mathbb{N}}
P\left(R_{\alpha'}(g)-R_{n,\alpha'}(g) \geq
\sqrt{\frac{\ln(n)+z+x_{m'}}{2n(1-\alpha')^2}}
+\frac{2}{1-\alpha'}\sqrt{\frac{V_{\mathcal{G}_{m'}}\ln(n+1)+\ln(2)}{n}}\right)\\
&\leq\displaystyle\sum_{\alpha'\in}\sum_{m'\in\mathbb{N}}
\frac{1}{n}e^{-z-x_{m'}}\leq \displaystyle\sum_{m'\in\mathbb{N}}
e^{-z-x_{m'}}\leq \Sigma e^{-z}.
\end{align*}

Considering now the complete interval, if
$\alpha'\in[0,\alpha_{max}]$, then $\exists\alpha''\in A$ such that
$\alpha''\leq\alpha'\leq\alpha''+\frac{1}{n}$. So from
\eqref{EqIniCorolarioCont}, with probability greater than
$\frac{1}{n}e^{-z-x_{m'}}$,
\begin{equation*}
\displaystyle\sup_{g\in\mathcal{G}_{m'}}\left(R_{\alpha''}(g)-R_{n,\alpha''}(g)\right)\leq
\sqrt{\frac{\ln(n)+z+x_{m'}}{2n(1-\alpha'')^2}}
+\frac{2}{1-\alpha''}\sqrt{\frac{V_{\mathcal{G}_{m'}}\ln(n+1)+\ln(2)}{n}},
\end{equation*}
then for all $\alpha'\in[0,\alpha_{max}]$
\begin{align*}
&\displaystyle\sup_{g\in\mathcal{G}_{m'}}\left(R_{\alpha'}(g)-R_{n,\alpha'}(g)\right)\\
&=\displaystyle\sup_{g\in\mathcal{G}_{m'}}\left(R_{\alpha'}(g)-R_{n,\alpha'}(g)+R_{\alpha''}(g)-R_{\alpha''}(g)+R_{n,\alpha''}(g)-R_{n,\alpha''}(g)\right)\\
&=\displaystyle\sup_{g\in\mathcal{G}_{m'}}\left([R_{\alpha''}(g)-R_{n,\alpha''}(g)]+[R_{\alpha'}(g)-R_{\alpha''}(g)]+[R_{n,\alpha''}(g)-R_{n,\alpha'}(g)]\right)\\
&\leq \sqrt{\frac{\ln(n)+z+x_{m'}}{2n(1-\alpha'')^2}}
+\frac{2}{1-\alpha''}\sqrt{\frac{V_{\mathcal{G}_{m'}}\ln(n+1)+\ln(2)}{n}}+\frac{1}{n(1-\alpha_{max})^2}\\
&\leq
\sqrt{\frac{\ln(n)+x_{m'}}{2n(1-\alpha')^2}}+\sqrt{\frac{z}{2n(1-\alpha')^2}}
+\frac{2}{1-\alpha'}\sqrt{\frac{V_{\mathcal{G}_{m'}}\ln(n+1)+\ln(2)}{n}}+\frac{1}{n(1-\alpha_{max})^2}.
\end{align*}
Where the next-to-last inequality comes from applying Proposition
\ref{DifEntreErroresRecortados} and that $R_{n,\alpha'}(g)\leq
R_{n,\alpha''}(g)$ and hence $R_{n,\alpha'}(g)-
R_{n,\alpha''}(g)\leq 0$ and the last one comes from
$\alpha''\leq\alpha'$. We can conclude that
\begin{equation*}
R_{\hat{\alpha}}(\mathcal{G}_{\hat{m}})-R_{n,\hat{\alpha}}(\mathcal{G}_{\hat{m}})\leq
\sqrt{\frac{\ln(n)+z+x_{\hat{m}}}{2n(1-\hat{\alpha})^2}}
+\frac{2}{1-\hat{\alpha}}\sqrt{\frac{V_{\mathcal{G}_{\hat{m}}}\ln(n+1)+\ln(2)}{n}}+\frac{1}{n(1-\alpha_{max})^2}.
\end{equation*}

\noindent Going back to \eqref{EqBasicaTeoCota3FuncFin}, except in a
set of probability not greater than $(\Sigma+1)e^{-z}$
\begin{eqnarray*}
R_{\hat{\alpha}}(\mathcal{G}_{\hat{m}}) &\leq&
R_{\alpha}(\mathcal{G}_m) + pen(\alpha,\mathcal{G}_m) +
\frac{\sqrt{R(\mathcal{G}_m)}}{\sqrt{2n}(1-\alpha)}
+\sqrt{\frac{z}{2n(1-\alpha)^2}} -
pen(\hat{\alpha},\mathcal{G}_{\hat{m}})\\
&+& \sqrt{\frac{\ln(n)+x_{\hat{m}}}{2n(1-\hat{\alpha})^2}}
+\sqrt{\frac{z}{2n(1-\hat{\alpha})^2}}
+\frac{2}{1-\hat{\alpha}}\sqrt{\frac{V_{\mathcal{G}_{\hat{m}}}\ln(n+1)+\ln(2)}{n}}+\frac{1}{n(1-\alpha_{max})^2}.
\end{eqnarray*}
Considering
\begin{equation*}
pen(\alpha,\mathcal{G}_m)=\sqrt{\frac{\ln(n+1)+x_m}{2n(1-\alpha)^2}}+
\frac{1}{(1-\alpha)}\sqrt{\frac{V_{\mathcal{G}_m}\ln(n+1)+ln(2)}{n}},
\end{equation*}
we have
\begin{eqnarray*}
R_{\hat{\alpha}}(\mathcal{G}_{\hat{m}}) &\leq&
R_{\alpha}(\mathcal{G}_m) + pen(\alpha,\mathcal{G}_m) +
\frac{\sqrt{R(\mathcal{G}_m)}}{\sqrt{2n}(1-\alpha)}
+\sqrt{\frac{z}{2n(1-\alpha)^2}}
+\sqrt{\frac{z}{2n(1-\hat{\alpha})^2}}+\frac{1}{n(1-\alpha_{max})^2}\\
&\leq&R_{\alpha}(\mathcal{G}_m) + pen(\alpha,\mathcal{G}_m) +
\frac{\sqrt{R(\mathcal{G}_m)}}{\sqrt{2n}(1-\alpha)}
+\sqrt{\frac{2z}{n(1-\frac{k_0}{n})^2}}+\frac{1}{n(1-\alpha_{max})^2}.
\end{eqnarray*}
Now grouping and integrating with respect to z,
\begin{eqnarray*}
E(R_{\hat{\alpha}}(\mathcal{G}_{\hat{m}}))&\leq& \min_{(\alpha,m)
\in [0,\alpha_{max}]\times\mathbb{N}}
\left(R_\alpha(\mathcal{G}_m)+pen(\alpha,\mathcal{G}_m)+
\frac{\sqrt{R(\mathcal{G}_m)}}{\sqrt{2n}(1-\alpha)}\right)\\
&+&\frac{1+\Sigma}{2(1-\alpha_{max})}\sqrt{\frac{\pi}{2n}}+
\frac{1}{n(1-\alpha_{max})^2}.
\end{eqnarray*}
\hfill $\Box$



\end{document}